# Effective Methods of QR-Decompositions of Square Complex Matrices by Fast Discrete Signal-Induced Heap Transforms

Artyom M. Grigoryan

[1] Department of Electrical and Computer Engineering,
The University of Texas at San Antonio, USA
amgrigoryan@utsa.edu

## *ABSTRACT*

The purpose of this work is to present an effective tool for computing different QR-decompositions of a complex nonsingular square matrix. The concept of the discrete signal-induced heap transform (DsiHT, Grigoryan 2006) is used. This transform is fast, has a unique algorithm for any length of the input vector/signal and can be used with different complex basic 2×2 transforms. The DsiHT is zeroing all components of the input signal while moving or heaping the energy of the signal to one component, for instance the first one. We describe three different types of QR-decompositions that use the basic transforms with the $\boldsymbol{T}$, $\boldsymbol{G}$, and $\boldsymbol{M}$-type complex matrices we introduce, as well as without matrices but using analytical formulas. We also present the mixed QR-decomposition, when different type DsiHTs are used in different stages of the algorithm. The number of such decompositions is greater than $3^{(N-1)}$, for an $N \times N$ complex matrix. Examples of the QR-decomposition are described in detail for the 4×4 and 6×6 complex matrices and compared with the known method of Householder transforms. The precision of the QR-decompositions of $N \times N$ matrices, when $N$ are 6, 13, 17, 19, 21, 40, 64, 100, 128, 201, 256, and 400 is also compared. The MATLAB-based scripts of the codes for QR-decompositions by the described DsiHTs are given.

## 1. INTRODUCTION

The QR-decomposition, or factorization of a non-singular matrix $\boldsymbol{A} = \boldsymbol{QR}$ into a unitary matrix $\boldsymbol{Q}$ and an upper triangular matrix $\boldsymbol{R}$, as well as the factorization $\boldsymbol{A} = \boldsymbol{QL}$ with a low triangular matrix $\boldsymbol{L}$ are powerful tools for solving linear systems of equations $\boldsymbol{y} = \boldsymbol{Ax}$ in many applications in computing and data analysis [1]-[7]. Here, the matrix $\boldsymbol{A}$ is a real or complex non-singular matrix. The matrix $\boldsymbol{Q}$ is unitary, and therefore its inverse is the transpose conjugate $\boldsymbol{Q}^*$ matrix. The calculation of inverse of the triangular matrix is not a difficult task. Therefore, the solution of the system of equation can be calculated by $\boldsymbol{x} = \boldsymbol{A}^{-1}\boldsymbol{y}$=$\boldsymbol{R}^{-1}\boldsymbol{Q}^*\boldsymbol{y}$ for the QR-decomposition, and $\boldsymbol{x} = \boldsymbol{A}^{-1}\boldsymbol{y}$=$\boldsymbol{L}^{-1}\boldsymbol{Q}^*\boldsymbol{y}$ for the QL-decomposition. Many known methods of QR-decomposition of real matrices were modified for the complex case. They include the Gramm-Schmidt process [8], the method of Householder transformations (or Householder reflections) [9], and the Givens rotations [10,11].

In this work we focus on the QR-decomposition and describe three types of decompositions, by using the concept of the discrete signal-induced heap transform (DsiHT) [12]-[14]. In the case of real matrices, the detail description of the DsiHT method of QR-decomposition is given in [16]. For complex matrices, there are different types of 2×2 basic complex unitary transforms can, that transfer energy of the 2-point signal to the first components, while zeroing the second one. The path of the DsiHT, i.e., the order of sequential processing (or rotating in some cases) of data of the signal is an important characteristic of the transform. The DsiHTs with different paths result in different QR-decomposition of the same matrix, as shown in [17] on examples with the so-called weak and strong-DsiHTs. The interesting property of the QR-decomposition by the DsiHT is in the presence of analytical equations that allow to calculate the transforms and their matrices without using the basic matrices of rotations. The case with complex matrices is much richer than real matrices, since there are many different basic transforms, not only "Givens rotations," that can be considered in the DsiHT. Examples of such transforms and their application in QR-decomposition of complex matrices are described and compared with the complex Householder transform-based QR-decomposition.

The rest of this paper is organized in the following way. In Section II, the concept of the DsiHT is described with examples of two-wheel carriages that illustrate the performance of the transform. The basic complex matrices of the 2×2 transforms composing the $N$-point DsiHT are described in Section III. These matrices are classified by the basic transforms use them and are called the $\boldsymbol{T}$, $\boldsymbol{G}$, and $\boldsymbol{M}$-type complex matrices. The concept of the DsiHT which can be calculated without matrices but by analytical formulas is also described. The QR-decompositions with three types of DsiHTs are presented with examples in Section IV. It is shown that the $\boldsymbol{M}$-type DsiHT based QR decomposition is close to the result of the Householder transform method, and other two types of QR-decomposition differ much. The scripts for MATLAB-based codes for computing the DsiHT and described QR-decompositions are given in Section V. In Section VI, the concept of the mixed QR-decomposition is presented, when different type DsiHTs are used in different stages of decomposition. The number of such decompositions is greater than $3^{(N-1)}$ for an $N \times N$ complex matrices, which is a very large number for large $N$. The questions related to the selection of the QR-decomposition from such large number of cases are not discussed here in detail, since it beyond the score of this work.

## 2. Basic 2×2 matrices in complex algebra

In this section, we describe briefly the concept of the discrete signal-induced heap transform (DsiHT) [12,15]. The transform is unitary and is defined by a non-zero vector, or signal $\boldsymbol{x} = (x_0, x_1, x_2, \dots, x_{N-1})$, without any constrain on the length $N$ and signal itself; it may be real and complex. This signal is called *the generator* of the DsiHT; the signal generates the unitary transform which is applying on other signals $\boldsymbol{z} = (z_0, z_1, z_2, \dots, z_{N-1})$. We consider the case when the $N$-point DsiHT is calculating by $(N-1)$ basic transformations $T$, each of which is applying only on two components of the renewal vector $\boldsymbol{z}$ in a certain order, or a path.

In the simple form, the DsiHT is calculated by applying a set of basic transformations $T$. The $2 \times 2$ matrix of such a transformation is defined by a chosen vector $(x, y)$ from the condition

$$\boldsymbol{T}\begin{bmatrix} x \\ y \end{bmatrix} = \begin{bmatrix} x' \\ 0 \end{bmatrix} \quad \text{and} \quad |x'|^2 = \sqrt{x^2 + y^2}. \tag{1}$$

In the case when $x$ and $y$ are real, $\boldsymbol{T}$ can be considered as the Givens rotation with the matrix

$$\boldsymbol{T} = \boldsymbol{T}_\varphi = \begin{bmatrix} \cos\varphi & -\sin\varphi \\ \sin\varphi & \cos\varphi \end{bmatrix}, \qquad \varphi = -\arctan\left(\frac{y}{x}\right).$$

If $x = 0$, the angle of rotation $\varphi = -\pi/2$, or $\pi/2$.

The concept of the $N$-point DsiHT of the signal $\boldsymbol{z}$ is illustrated in the diagram of Fig. 1. This is the so-called weak carriage with two wheels, one "rotates" the generator $\boldsymbol{x}$ and another wheel "rotates" the signal $\boldsymbol{z}$. When the carriage moves, the generator components work together with the first element $x_0$ and update/renew its value at each step. At the same time, the transformation $T$, which is determined during the rotation of the first wheel, is applying to the two components of the input signal $\boldsymbol{z}$. In this wheel, the components of $\boldsymbol{z}$ are processing together with the first element $z_0$ and update its value.

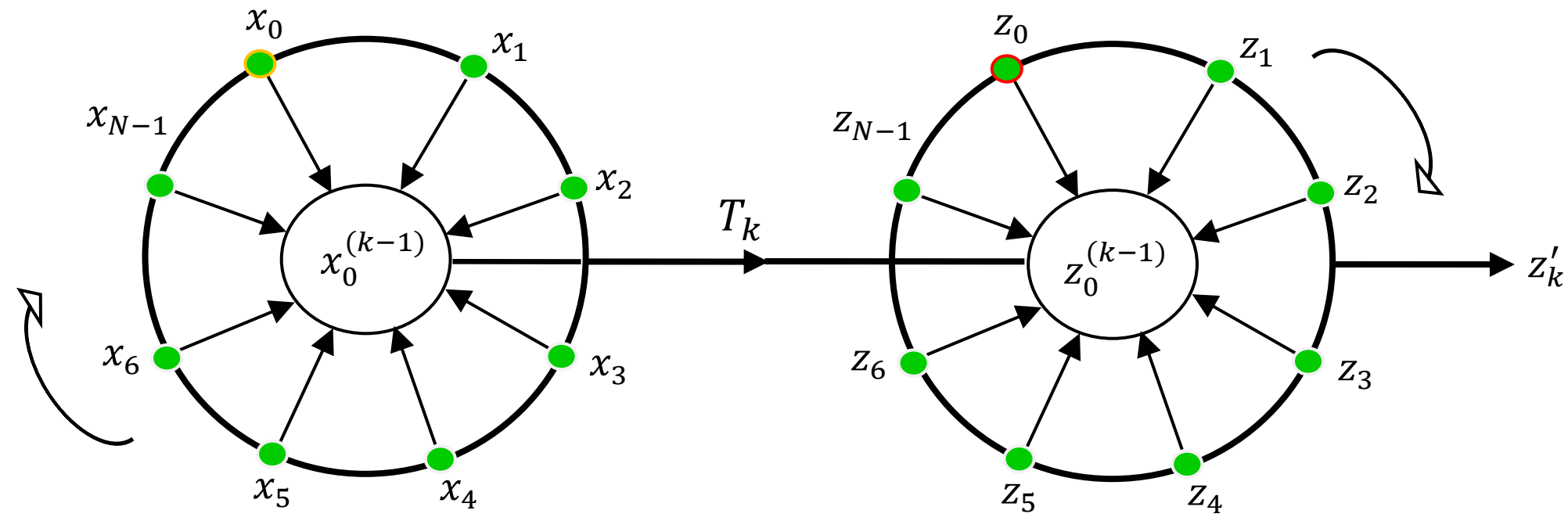

Figure 1: The two-wheel carriage of the DsiHT.

The $N$-point DsiHT with such a carriage is called a weak DsiHT, since the components $x_n$ and $z_n$ of the generator and signal are processed in the natural order of the index $n$, i.e., $x_0$ with $x_1$, then $x_2, x_3, \dots, x_{N-1}$, and the same for the signal $\boldsymbol{z}$. This is the natural path of the DsiHT and the path can be chosen differently [16]. For instance, the concept of the strong DsiHT is defined by the path in order $x_{N-1}$ with $x_{N-2}$, then $x_{N-3}, \dots, x_1, x_0$, and the same for the signal $\boldsymbol{z}$. The carriage of the strong DsiHT is illustrated in Fig. 2.

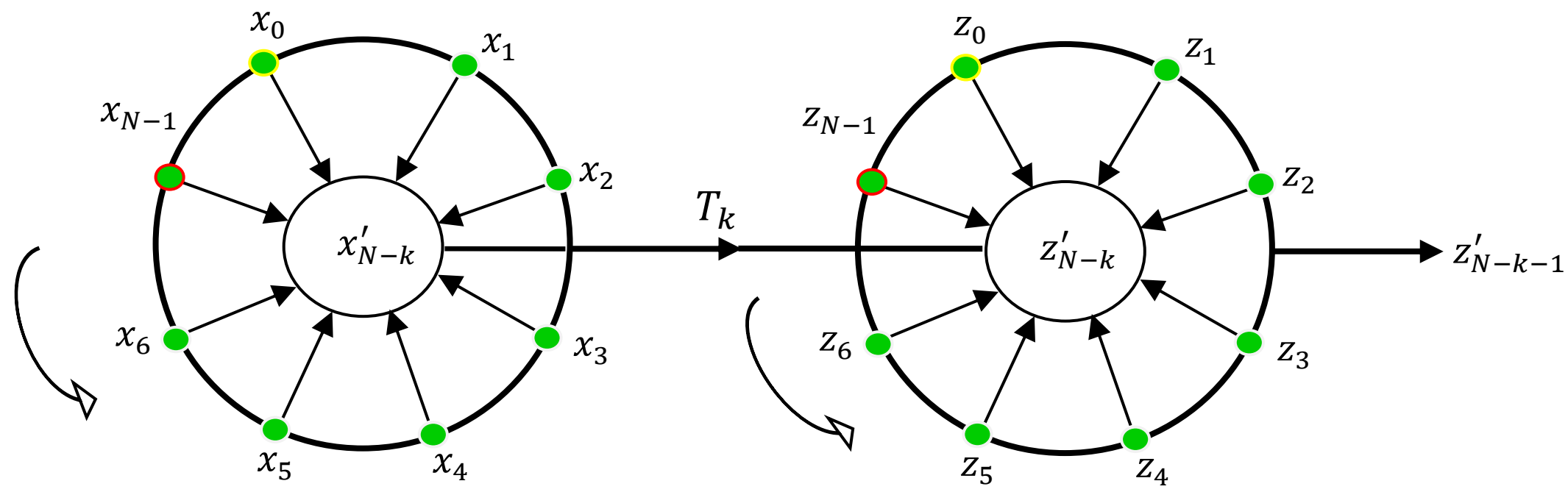

Figure 2: The two-wheel carriage of the strong DsiHT.

We consider the DsiHT with the natural path.

**Algorithm of the DsiHT** with the generator $\boldsymbol{x} = (x_0, x_1, x_2, \dots,, x_{N-1})$.

The input signal is $\boldsymbol{z} = (z_0, z_1, z_2, \dots,, z_{N-1})$.

1. Stage $k = 1$.
   - Calculate the matrix of the transform $\boldsymbol{T}_1$ that is generated by the vector $(x_0, x_1)'$.
   - Calculate the transform $\boldsymbol{T}_1: (z_0, z_1)' \to \left(z_0^{(1)}, z_1'\right)'$.
   - Calculate the new value $x_0^{(1)}$ of $x_0$ in the transform $\boldsymbol{T}_1: (x_0, x_1)' \to \left(x_0^{(1)}, 0\right)'$.
2. Stage $k = 2$.
   - Calculate the matrix of the transform $\boldsymbol{T}_k$ that is generated by the vector $(x_0^{(k-1)}, x_k)'$.
   - Calculate the transform $\boldsymbol{T}_k: \left(z_0^{(k-1)}, z_k\right)' \to \left(z_0^{(k)}, z_k'\right)'$.
   - Calculate the new value $x_0^{(k)}$ in the transform $\boldsymbol{T}_2: \left(x_0^{(k-1)}, x_k\right)' \to \left(x_0^{(k)}, 0\right)'$.
3. Stage $k = 3,4, \dots, N-1$.
   - Continue calculations as in Step 2, to get the rest values of the output $z_3', z_4', \dots, z_{N-1}'$, and $z_0^{(N)}$.

The output signal is $\boldsymbol{z}' = (z_0^{(N)}, z_1', z_2', z_3', \dots,, z_{N-1}')$.

To simplicity the notations in the two-wheel carriage in Fig. 1, $x_0'$ denotes $x_0^{(k)}$ and $z_0'$ denotes $z_0^{(k)}$ on the kth stage of rotations. Thus, during the composition of the $N$-point DsiHT, $H$, a set of parameters, or angles $\varphi_1, \varphi_2, \dots, \varphi_{N-1}$, is calculated, which is called *the angular representation* of the generator $\boldsymbol{x}$ [15]. The transformation $H$ is separable and calculated as

$$H = T_{\varphi_{N-1}} \dots T_{\varphi_3} T_{\varphi_2} T_{\varphi_1}. \tag{2}$$

When $\boldsymbol{z} = \boldsymbol{x}$, the transform of the vector $\boldsymbol{x}$ is collected to one heap; it is transferred to the first component, $T(\boldsymbol{x}) = \left(x_0^{(N-1)}, 0, 0, \dots, 0\right)$, where $x_0^{(N-1)} = \|\boldsymbol{x}\| = \sqrt{x_0^2 + x_1^2 + \dots + x_{N-1}^2}$.

***Example 1:*** For the $N = 6$ case, we consider the generator $\boldsymbol{x} = (1,1,2,4,3,1)$ with the norm $\|\boldsymbol{x}\| = \sqrt{32} = 5.6569$. The matrix and the set of five angles (in radians) of the 6-point $\boldsymbol{x}$-DsiHT are

$$\boldsymbol{H}_1 = \begin{pmatrix} 0.1768 & 0.1768 & 0.3536 & 0.7071 & 0.5303 & 0.1768 \\ -0.7071 & 0.7071 & 0 & 0 & 0 & 0 \\ -0.5774 & -0.5774 & 0.5774 & 0 & 0 & 0 \\ -0.3482 & -0.3482 & -0.6963 & 0.5222 & 0 & 0 \\ -0.1149 & -0.1149 & -0.2298 & -0.4595 & 0.8424 & 0 \\ -0.0318 & -0.0318 & -0.0635 & -0.1270 & -0.0953 & 0.9843 \end{pmatrix},$$

$$\{\varphi_1,\varphi_2,\varphi_3,\varphi_4,\varphi_5\} = \{-0.7854, -0.9553, -1.0213, -0.5690, -0.1777\}.$$

The matrix $\boldsymbol{H}_1$ of the transform can be written as $\boldsymbol{H}_1 = \boldsymbol{D}_1\boldsymbol{A}_1$, where the diagonal matrix

$$\boldsymbol{D}_1 = \text{diag}\{0.1768, -0.7071, -0.5774, -0.3482, -0.1149, -0.0318\}$$

and the matrix $\boldsymbol{A}_1$ is

$$\boldsymbol{A}_1 = \begin{pmatrix} 1 & 1 & 2 & 4 & 3 & 1 \\ 1 & -1 & 0 & 0 & 0 & 0 \\ 1 & 1 & -1 & 0 & 0 & 0 \\ 1 & 1 & 2 & -1.5 & 0 & 0 \\ 1 & 1 & 2 & 4 & -22/3 & 0 \\ 1 & 1 & 2 & 4 & 3 & -31 \end{pmatrix}.$$

The determinant of the matrix $\det \boldsymbol{H}_1 = 1$. The first row of the matrix is the generator $\boldsymbol{x}$ and the last row is the same generator only with the splash at the end. This value is calculated by $31 = 1 + 1 + 2^2 + 4^2 + 3^2$; it is the energy of the generator without the last component, i.e., $\|\boldsymbol{x}\|^2 - 1$.

For the 6-point DsiHT with the path of the strong carriage, we have the following data. The angular representation is

$$\{\varphi_1,\varphi_1,\dots,\varphi_6\} = \{-1.3931, -1.3902, -1.1970, -0.6690, -0.3218\},$$

and the matrix $\boldsymbol{H}_2$ of this DsiHT with $\det \boldsymbol{H}_2 = 1$ is

$$\boldsymbol{H}_2 = \begin{pmatrix} 0.1768 & 0.1768 & 0.3536 & 0.7071 & 0.5303 & 0.1768 \\ -0.9843 & 0.0318 & 0.0635 & 0.1270 & 0.0953 & 0.0318 \\ 0 & -0.9837 & 0.0656 & 0.1312 & 0.0984 & 0.0328 \\ 0 & 0 & -0.9309 & 0.2864 & 0.2148 & 0.0716 \\ 0 & 0 & 0 & -0.6202 & 0.7442 & 0.2481 \\ 0 & 0 & 0 & 0 & -0.3162 & 0.9487 \end{pmatrix}.$$

This matrix can be presented as

$$\boldsymbol{H}_2 = \boldsymbol{D}_2\boldsymbol{A}_2, \qquad \boldsymbol{A}_2 = \begin{pmatrix} 1 & 1 & 2 & 4 & 3 & 1 \\ -31 & 1 & 2 & 4 & 3 & 1 \\ 0 & -30 & 2 & 4 & 3 & 1 \\ 0 & 0 & -13 & 4 & 3 & 1 \\ 0 & 0 & 0 & -2.5 & 3 & 1 \\ 0 & 0 & 0 & 0 & -1/3 & 1 \end{pmatrix}.$$

Here, the diagonal matrix $\boldsymbol{D}_2 = \text{diag}\{0.1768, 0.0318, 0.0328, 0.0716, 0.2481, 0.9487\}$. Thus, we have two different matrices $\boldsymbol{H}_1$ and $\boldsymbol{H}_2$ such that

$$\boldsymbol{H}_1\boldsymbol{x}' = \boldsymbol{H}_2\boldsymbol{x}' = (\|\boldsymbol{x}\|, 0,0,0,0,0)' = (5.6569,0,0,0,0,0)'.$$

If we take a vector $\boldsymbol{z}$, for instance equal $\boldsymbol{z} = (4, -2, 3, -1, 7, 2)$ with the norm $\|\boldsymbol{z}\| = 9.1104$, we obtain the following transforms:

$$\boldsymbol{z}_1 = \boldsymbol{H}_1\boldsymbol{z}' = (4.7730, -4.2426, 0.5774, -3.3075, 5.4375, 1.1748)',$$
$$\boldsymbol{z}_2 = \boldsymbol{H}_2\boldsymbol{z}' = (4.7730, -3.2068, 2.7873, -1.4322, 6.3258, -0.3162)',$$

and $\|\boldsymbol{z}_1\| = \|\boldsymbol{z}_2\| = 9.1104 = \|\boldsymbol{z}\|$.

Below are the script of the MATLAB-based code to calculate the matrices of these two DsiHTs.

```
% ---------------------------------------------------------------
 x=[1 1 2 4 3 1]';
 [H1,U1]=matrix_msob(x);
 [H2,U2]=matrix_msob2(x);
```

```
% call: matrix_msob.m   / from the library of Art Grigoryan
% To calculate the matrix A and angles U of the X-generated DsiHT.
function [A,U]=matrix_msob(y)
    N=length(y);
    A=eye(N);
    pi2=pi/2;
    U=zeros(1,N);
    for i1=2:N
        yi1=y(i1);   y1=y(1);
        if abs(yi1) >eps
            if y1==0
                u=pi2;
            else
                u=acot(-y1/yi1);
            end
            if y1<=0 u=u+pi; end
            U(i1)=u;
            c1=cos(u); s1=sin(u);
            y(1)=[c1 -s1]*[y1 yi1]';
            S=eye(N);
            S(1,1) =c1;   S(i1,1)=-s1;
            S(1,i1)=s1;   S(i1,i1)=c1;
            A=A*S;
        end
    end
    U(1)=y(1);  A=A';
end
% ----------------------------------------------------------------------
% call: matrix_msob2.m   / from the library of Art Grigoryan
% To calculate the matrix A and angles U of the X-generated strong DsiHT.
function [A,U]=matrix_msob2(y)
    N=length(y);
    A=eye(N);    pi2=pi/2;
    U=zeros(1,N);
    for i1=N:-1:2
        y1=y(i1-1);  yi1=y(i1);
        if abs(yi1) >eps
            if y1==0
                u=pi2;
            else
                u=acot(-y1/yi1);
            end
            U(i1)=u;
            c1=cos(u); s1=sin(u);
            y(i1-1)=[c1 -s1]*[y1 yi1]';
            S=zeros(N);
            for k=1:N
                S(k,k)=1;
            end
            S(i1-1,i1-1) =c1;   S(i1,i1-1)=-s1;
            S(i1-1,i1)=s1;      S(i1,i1)=c1;
            A=A*S;
        end
    end
    U(1)=y(1);  A=A';
% --------------------------------------------------------------------
```

## 2.1 Calculation of the DsiHT by Analytical Equations

It is important to note that the $N$-point DsiHT can be obtained without calculation of the angles $\varphi_k$ and trigonometric functions $\cos\varphi_k$ and $\sin\varphi_k$, but by using the analytical formulas [15]. For that, we consider the following notations, which represent respectively the partial cross-correlation of $\boldsymbol{z}$ with the vector generator $\boldsymbol{x}$:

$$E_k(\boldsymbol{z},\boldsymbol{x}) = z_0x_0 + z_1x_1 + \cdots + z_kx_k, \tag{3}$$

and the partial and full energies of the signal generator

$$E_k^2(\boldsymbol{x}) = E_k(\boldsymbol{x},\boldsymbol{x}) = x_0^2 + x_1^2 + \cdots + x_k^2, \quad k = 0,1,\ldots,(N-1). \tag{4}$$

The components of the DsiHT on the $k$-th iteration can be expressed by the correlation data as follows:

$$z_k^{(1)} = \frac{E_k(\boldsymbol{x},\boldsymbol{x})z_k - E_k(\boldsymbol{z},\boldsymbol{x})x_k}{E_{k-1}(\boldsymbol{x})E_k(\boldsymbol{x})}, \qquad k = 1,2,\ldots,(N-1). \tag{5}$$

On the final stage, the value of the first component is defined by

$$z_0^{(N-1)} = E_{N-1}(\boldsymbol{z},\boldsymbol{x})/E_{N-1}(\boldsymbol{x}), \tag{6}$$

which is the correlation coefficient of the input signal $\boldsymbol{z}$ with the normalized generator $\boldsymbol{x}$. For a given generator, all values of $E_k(\boldsymbol{x},\boldsymbol{x})$ and $E_{k-1}(\boldsymbol{x})E_k(\boldsymbol{x})$, can be calculated in advance. In the case when $\boldsymbol{z}=\boldsymbol{x}$,

$$x_0^{(N-1)} = \frac{E_{N-1}(\boldsymbol{x},\boldsymbol{x})}{E_{N-1}^2(\boldsymbol{x})} = \|\boldsymbol{x}\|$$

and $x_k^{(1)} = 0$, for all $k = 1,2,\ldots,(N-1)$.

The coefficients $h_{n,m}$ of the matrix of the $N$-point DsiHT can be obtained from equations (3)-(6). The m-th column of the DsiHT matrix $\boldsymbol{H}$ can be calculated, by applying the unit vector $\boldsymbol{e}_m = (0,0,\ldots,1,\ldots0)'$ with 1 on the m-th position, where $m \in \{0,1,\ldots,N-1\}$. Therefore, the coefficients of the transform can be calculated by

$$h_{0,m} = \frac{E_{N-1}(\boldsymbol{e}_m,\boldsymbol{x})}{E_{N-1}(\boldsymbol{x})}, \quad m = 0:(N-1), \tag{7}$$

$$h_{n,m} = \frac{E_n^2(\boldsymbol{x})(\boldsymbol{e}_m)_n - E_n(\boldsymbol{e}_m,x)x_n}{E_{n-1}(\boldsymbol{x})E_n(\boldsymbol{x})}, \tag{8}$$

where $n = 1:(N-1)$.

## 3. COMPLEX BASIC MATRICES

The $N$-point DsiHT is calculated, by applying $(N-1)$ basic transformations $T$ on two different components of the renewal vector in a certain order, or the path. In this section, we consider the concept of the complex DsiHT. The basic transformation of the DsiHT, which is defined by a complex vector $(x_0,x_1)'$, and then, is applied to a complex input $(z_0,z_1)'$ is calculated as follows:

$$T: \begin{bmatrix} z_0 \\ z_1 \end{bmatrix} \to \begin{bmatrix} z_0' \\ z_1' \end{bmatrix} = \frac{\text{sign}(\text{Real}(x_0))}{\sqrt{|x_0|^2 + |x_1|^2}} \begin{bmatrix} \bar{x}_0 & \bar{x}_1 \\ -x_1 & x_0 \end{bmatrix} \begin{bmatrix} z_0 \\ z_1 \end{bmatrix}. \tag{9}$$

The complex matrix of this transformation is

$$\boldsymbol{T} = \frac{\text{sign}(\text{Real}(x_0))}{\sqrt{|x_0|^2 + |x_1|^2}} \begin{bmatrix} \bar{x}_0 & \bar{x}_1 \\ -x_1 & x_0 \end{bmatrix}$$

and the determinant is 1. It is not difficult to verify that the matrix $\boldsymbol{T}$ is unitary, i.e., the multiplication of $\boldsymbol{T}$ with its conjugate transpositon $\boldsymbol{T}(\overline{\boldsymbol{T}}^T)= \boldsymbol{I}$, where $\boldsymbol{I}$ is the 2×2 identity matrix.

When the input equals to the generator, i.e., $(z_0,z_1)' = (x_0,x_1)'$, we obtain the real transform

$$T: \begin{bmatrix} x_0 \\ x_1 \end{bmatrix} \to \text{sign}(\text{Real}(x_0)) \begin{bmatrix} \sqrt{|x_0|^2 + |x_1|^2} \\ 0 \end{bmatrix}. \tag{10}$$

## 3.1. DsiHT with Analytical Equations

The complex DsiHT can also be calculated by using analytical equations (5) and (6), similar to the case with real vectors. For complex vectors, the partial cross-correlation of $\boldsymbol{z}$ with the vector- generator $\boldsymbol{x}$ in the first $k$ components is calculated by

$$E_k(\boldsymbol{z},\boldsymbol{x}) = z_0\bar{x}_0 + z_1\bar{x}_1 + \cdots + z_k\bar{x}_k, \tag{11}$$

and the energies in the first $k$ components of the signal-generator by

$$E_k^2(\boldsymbol{x}) = E_k(\boldsymbol{x},\boldsymbol{x}) = |x_0|^2 + |x_1|^2 + \cdots + |x_k|^2, \quad k = 0,1,\dots,(N-1). \tag{12}$$

The 2×2 matrix calculated by analytical equations (7) and (8) for the $N = 2$ case will be denoted by $\boldsymbol{M}$. This matrix is different from matrix $\boldsymbol{T}$ and equals

$$\boldsymbol{M} = \frac{1}{\sqrt{|x_0|^2 + |x_1|^2}} \begin{bmatrix} \bar{x}_0 & \bar{x}_1 \\ -x_1\dfrac{\bar{x}_0}{|x_0|} & |x_0| \end{bmatrix}. \tag{13}$$

One coefficient of the matrix is a real number, and all coefficients are complex in the matrix $\boldsymbol{T}$. The matrix $\boldsymbol{M}$ is unitary and the determinant of the matrix is the complex number

$$\det\boldsymbol{M} = \frac{1}{|x_0|^2 + |x_1|^2}\left[\bar{x}_0|x_0| + \bar{x}_1 x_1\frac{\bar{x}_0}{|x_0|}\right] = \frac{\bar{x}_0}{|x_0|} \tag{14}$$

and $|\det\boldsymbol{M}| = 1$. The matrix product $\boldsymbol{M}(x_0,x_1)'$ equals

$$\frac{1}{\sqrt{|x_0|^2 + |x_1|^2}} \begin{bmatrix} \bar{x}_0 & \bar{x}_1 \\ -x_1\dfrac{\bar{x}_0}{|x_0|} & |x_0| \end{bmatrix}\begin{bmatrix} x_0 \\ x_1 \end{bmatrix} = \frac{1}{\sqrt{|x_0|^2 + |x_1|^2}}\begin{bmatrix} |x_0|^2 + |x_1|^2 \\ 0 \end{bmatrix} = \begin{bmatrix} \sqrt{|x_0|^2 + |x_1|^2} \\ 0 \end{bmatrix}.$$

## 3.2. DsiHT with Complex Givens Rotations

We consider the known complex Givens rotation [4] which is defined by the matrix

$$\boldsymbol{G} = \frac{1}{\sqrt{|x_0|^2 + |x_1|^2}}\begin{bmatrix} |x_0| & \text{sign}(x_0)\bar{x}_1 \\ -\overline{\text{sign}}(x_0)x_1 & |x_0| \end{bmatrix} = \frac{1}{\sqrt{|x_0|^2 + |x_1|^2}}\begin{bmatrix} |x_0| & \dfrac{x_0}{|x_0|}\bar{x}_1 \\ -x_1\dfrac{\bar{x}_0}{|x_0|} & |x_0| \end{bmatrix}. \tag{15}$$

Here, the complex sign function is defined by $\text{sign}(x_0) = x_0/|x_0|$. The determinant of this matrix is 1. For the $x_0 = 0$ case, the above matrix is defined as

$$\boldsymbol{G} = \frac{1}{|x_1|}\begin{bmatrix} 0 & \bar{x}_1 \\ -x_1 & 0 \end{bmatrix},$$

i.e., it is considered that $\text{sign}(0) = 1$. The matrix $\boldsymbol{G}$ is half-complex, meaning that in each row and column there is a real number $|x_0|$. It should be noted that all coefficients of the basic matrix $\boldsymbol{T}$ in Eq. 9 are complex numbers. When applying the matrix $\boldsymbol{G}$ on the vector-generator $(x_0,x_1)'$ we obtain

$$\boldsymbol{G}\begin{bmatrix} x_0 \\ x_1 \end{bmatrix} = \frac{1}{\sqrt{|x_0|^2 + |x_1|^2}}\begin{bmatrix} |x_0|x_0 + \dfrac{x_0}{|x_0|}|x_1|^2 \\ 0 \end{bmatrix} = \frac{x_0}{|x_0|}\begin{bmatrix} \sqrt{|x_0|^2 + |x_1|^2} \\ 0 \end{bmatrix},$$

or

$$\boldsymbol{G}\begin{bmatrix} x_0 \\ x_1 \end{bmatrix} = \text{sign}(x_0)\begin{bmatrix} \sqrt{|x_0|^2 + |x_1|^2} \\ 0 \end{bmatrix}, \tag{16}$$

where $\text{sign}(x_0)$ is a complex number with norm 1. The matrix $\boldsymbol{G}$ is unitary, i.e., $\boldsymbol{G}(\overline{\boldsymbol{G}^T}) = \boldsymbol{I}$.

It should be noted that one can express the matrix $\boldsymbol{G}$ by $\boldsymbol{T}$ and by $\mathbf{M}$, as

$$\boldsymbol{G}=\operatorname{sign}\big(\operatorname{Real}(x_0)\big)\begin{bmatrix}\dfrac{x_0}{|x_0|} & 0\\ 0 & \dfrac{\bar{x}_0}{|x_0|}\end{bmatrix}\boldsymbol{T},\qquad \boldsymbol{G}=\begin{bmatrix}\dfrac{x_0}{|x_0|} & 0\\ 0 & 1\end{bmatrix}\boldsymbol{M}. \tag{17}$$

***Example 2:*** Let us consider two complex numbers $x_0=1+3i$ and $x_1=-2+5i$. For these numbers, $|x_0|^2=1+3^2=10$ and $|x_1|^2=2^2+5^2=29$. We obtain the following matrices:

$$\boldsymbol{T}=\frac{\operatorname{sign}(\operatorname{Real}(x_0))}{\sqrt{|x_0|^2+|x_1|^2}}\begin{bmatrix}\bar{x}_0 & \bar{x}_1\\ -x_1 & x_0\end{bmatrix}=\frac{1}{\sqrt{39}}\begin{bmatrix}1-3i & -2-5i\\ 2-5i, & 1+3i\end{bmatrix},$$

$$\boldsymbol{M}=\frac{1}{\sqrt{|x_0|^2+|x_1|^2}}\begin{bmatrix}\bar{x}_0 & \bar{x}_1\\ -x_1\dfrac{\bar{x}_0}{|x_0|} & |x_0|\end{bmatrix}=\frac{1}{\sqrt{39}}\begin{bmatrix}1-3i & -2-5i\\ \dfrac{-13-11i}{\sqrt{10}} & \sqrt{10}\end{bmatrix},$$

and

$$\boldsymbol{G}=\frac{1}{\sqrt{|x_0|^2+|x_1|^2}}\begin{bmatrix}|x_0| & \dfrac{x_0}{|x_0|}\bar{x}_1\\ -x_1\dfrac{\bar{x}_0}{|x_0|} & |x_0|\end{bmatrix}=\frac{1}{\sqrt{39}}\begin{bmatrix}\sqrt{10} & \dfrac{13-11i}{\sqrt{10}}\\ \dfrac{-13-11i}{\sqrt{10}} & \sqrt{10}\end{bmatrix}.$$

The determinants of these matrices, $\det\boldsymbol{T}=\det\boldsymbol{G}=1$, $\det\boldsymbol{M}=\ 0.3162-0.9487i$, and $|\det\boldsymbol{M}|=1$. Up to the coefficient $1/\sqrt{39}$, the matrix $\boldsymbol{T}$ is integer-valued for integer complex generator, and matrices $\boldsymbol{G}$ and $\boldsymbol{M}$ are more complicated; they have additional coefficents with $\sqrt{10}=3.1623$. When applying these matrices on the vector $(x_0,x_1)'=(1+3i,-2+5i)'$, we obtain the following vectors:

$$\boldsymbol{T}\begin{bmatrix}x_0\\ x_1\end{bmatrix}=\boldsymbol{M}\begin{bmatrix}x_0\\ x_1\end{bmatrix}=\begin{bmatrix}6.2450\\ 0\end{bmatrix},\quad \boldsymbol{G}\begin{bmatrix}x_0\\ x_1\end{bmatrix}=\begin{bmatrix}1.9748+5.9245i\\ 0\end{bmatrix},$$

and $|1.9748+5.9245i|=6.2450$. Thus, the results of the transforms $\boldsymbol{T}$ and $\boldsymbol{M}$ are different from $\boldsymbol{G}$ in the first component which is complex. Now, we consider the transforms of the vector $z=(z_0,z_1)'=(-7+2i,3-5i)'$,

$$\boldsymbol{T}\begin{bmatrix}z_0\\ z_1\end{bmatrix}=\begin{bmatrix}-5.1241+2.8823i\\ 2.2418\ +\ 6.8855i\end{bmatrix},\qquad \boldsymbol{M}\begin{bmatrix}z_0\\ z_1\end{bmatrix}=\begin{bmatrix}-5.1241+2.8823i\\ 7.2411+0.0506i\end{bmatrix},$$

$$\boldsymbol{G}\begin{bmatrix}z_0\\ z_1\end{bmatrix}=\begin{bmatrix}-4.3548-3.9497i\\ 7.2411+0.0506i\end{bmatrix}.$$

The outputs of $\boldsymbol{T}$, $\boldsymbol{M}$, and $\boldsymbol{G}$ are different. The matrix $\boldsymbol{M}$ has commom with matrices $\boldsymbol{T}$ and $\boldsymbol{G}$. The second components of the transform $\boldsymbol{G}$ and $\boldsymbol{M}$ are the same, and the first components of the transform $\boldsymbol{T}$ and $\boldsymbol{M}$ are the same. For all these transforms, the energy of the input vector is preserved,

$$\begin{gathered}5.1241^2+2.8823^2=4.3548^2+3.9497^2=34.5641,\\ 7.2411^2+0.0506^2=52.4359,\\ 2.24181^2\ +\ 6.88551^2=52.4359,\\ 34.5641+52.4359=87=(7^2+2^2)+(3^2+5^2)=|\mathrm{z}|^2.\end{gathered}$$

To describe the difference between the DsiHT generated by basic 2×2 transforms of different types, we consider the matrices of these DsiHTs for the $N=4$ case.

***Example 3:*** The complex vector-generator is

$$x=(x_0,x_1,x_2,x_3)=(7+4i,3+7\mathrm{i},-6+2\mathrm{i},1+2\mathrm{i}).$$

We have the following 4×4 matrices of the DsiHT calculated with basic transforms $\boldsymbol{T}$, $\boldsymbol{M}$, and $\boldsymbol{G}$:

$$\boldsymbol{T}_{4\times4}=\begin{bmatrix} 0.5401 & 0.2315 & -0.4629 & 0.0772\\ -0.2705 & 0.6312 & 0 & 0\\ 0.2401 & 0.0282 & 0.8687 & 0\\ -0.0906 & -0.1027 & 0.0121 & 0.9850 \end{bmatrix}+i\begin{bmatrix} -0.3086 & -0.5401 & -0.1543 & -0.1543\\ -0.6312 & 0.3607 & 0 & 0\\ -0.2684 & -0.3390 & 0 & 0\\ -0.0604 & 0.0060 & 0.0846 & 0 \end{bmatrix},$$

$$\boldsymbol{M}_{4\times4}=\begin{bmatrix} 0.5401 & 0.2315 & -0.4629 & 0.0772\\ -0.5480 & 0.7269 & 0 & 0\\ 0.2401 & 0.0282 & 0.8687 & 0\\ -0.0906 & -0.1027 & 0.0121 & 0.9850 \end{bmatrix}+i\begin{bmatrix} -0.3086 & -0.5401 & -0.1543 & -0.1543\\ -0.4138 & 0 & 0 & 0\\ -0.2684 & -0.3390 & 0 & 0\\ -0.0604 & 0.0060 & 0.0846 & 0 \end{bmatrix},$$

$$\boldsymbol{G}_{4\times4}=\begin{bmatrix} 0.6220 & 0.4687 & -0.3254 & 0.1435\\ -0.5480 & 0.7269 & 0 & 0\\ 0.2401 & 0.0282 & 0.8687 & 0\\ -0.0906 & -0.1027 & 0.0121 & 0.9850 \end{bmatrix}+i\begin{bmatrix} 0 & -0.3541 & -0.3636 & -0.0957\\ -0.4138 & 0 & 0 & 0\\ -0.2684 & -0.3390 & 0 & 0\\ -0.0604 & 0.0060 & 0.0846 & 0 \end{bmatrix}.$$

The determinants of these matrices equal $\det\boldsymbol{T}=\det\boldsymbol{G}=1$, $\det\boldsymbol{M}=\ 0.8622-0.4961i$, and $|\det\boldsymbol{M}|=1$. The matrices $\boldsymbol{T}_{4\times4}$ and $\boldsymbol{M}_{4\times4}$ are different only in coefficients of the 2[nd] rows. The difference of the matrices $\boldsymbol{T}_{4\times4}$ and $\boldsymbol{G}_{4\times4}$ is in the first two rows. In the matrices $\boldsymbol{T}_{4\times4}$ and $\boldsymbol{M}_{4\times4}$, the first rows are proportional to the vector generator. Indeed, for the matrix $\boldsymbol{T}_{4\times4}$ we have

$$(\boldsymbol{T}_{4\times4})_{0,0:3}=(\boldsymbol{M}_{4\times4})_{0,0:3}=\overline{(7+4i,3+7\mathrm{i},-6+2\mathrm{i},1+2\mathrm{i})}\,\frac{0.5401}{7}=(0.5401/7)\overline{\boldsymbol{x}}.$$

The vector-generator $\boldsymbol{x}$ is not in the first row of matrix $\boldsymbol{G}$. The first coefficient of this matrix is the real number 0.6220, not complex as in the vector-generator.

For the input vector $\boldsymbol{z}=(2-3i,1-4i,-7+i,3+5i)$, the transforms equal

$$\begin{aligned} \boldsymbol{T}_{4\times4}\boldsymbol{z}\ &=\ (2.6232-3.1632i,-0.3607-2.6148i,-7.7334-0.8404i,2.3447+4.9129i),\\ \boldsymbol{M}_{4\times4}\boldsymbol{z}&=(2.6232-3.1632,-1.6105-2.0914i,-7.7334-0.8404i,2.3447+4.9129i),\\ \boldsymbol{G}_{4\times4}\boldsymbol{z}&=(3.8469-1.4450i,-1.6105-2.0914i,-7.7334-0.8404i,2.3447+4.9129i). \end{aligned}$$

It is not difficult to verify that in the general $N\geq4$ case, the difference between the matrices $\boldsymbol{T}_{N\times N}$ and $\boldsymbol{G}_{N\times N}$ is only in the coefficeints of the first two rows. Thus, the DsiHT of signals differ only in the first components of the transform, if instead of the basic matrices $\boldsymbol{T}$ the matrices $\boldsymbol{G}$ are used.

As illustrative example, we consider the complex signal $\boldsymbol{z}$ of length 500 with real and imaginary parts shown in Fig. 3 in parts (a) and (b), respectively.

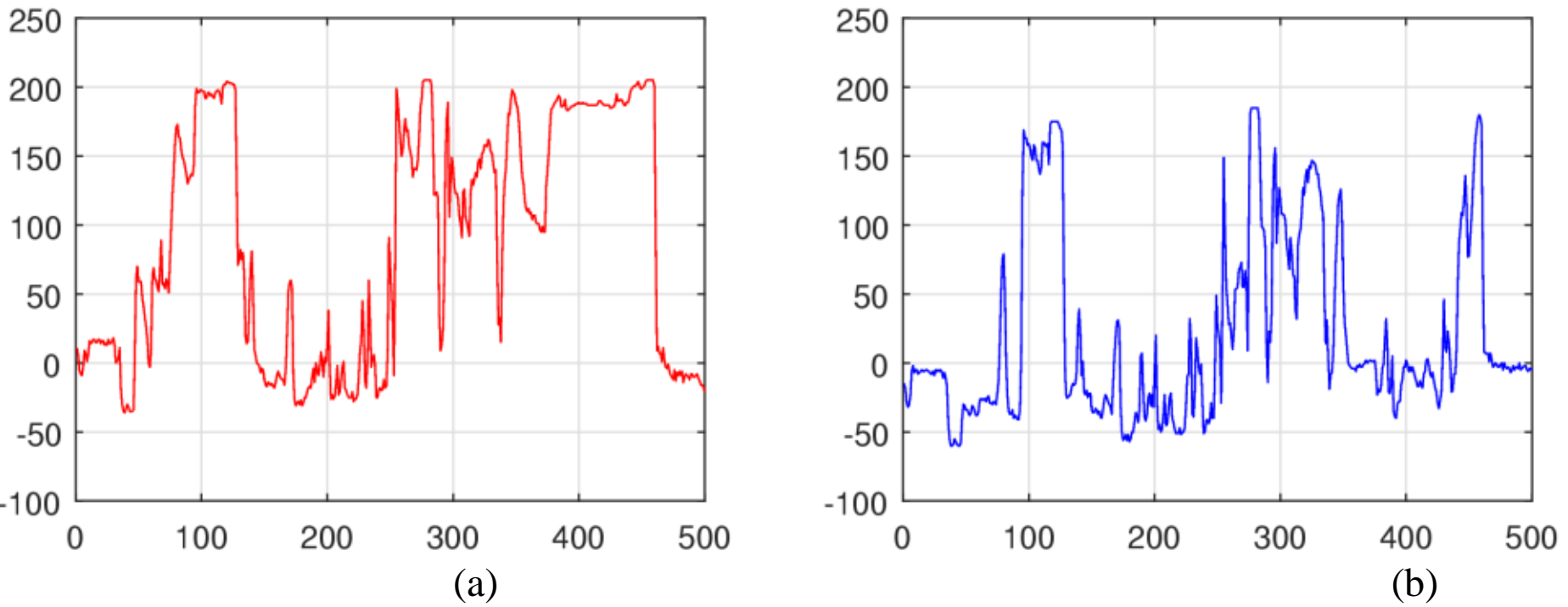

**Figure 3:** The complex signal $\boldsymbol{z}$ of length 500; (a) real and (b) imaginary parts.

The signal-generator $\boldsymbol{x}$ is calculated by

$$x(n)=4+i(-1)^2/8,\qquad n=0{:}\,499.$$

The complex signal $\boldsymbol{z}$ after processing by the DsiHT with generator $\boldsymbol{x}$ is shown in Fig. 4. The basic transforms of this DsiHT are with the matrices of type $\boldsymbol{T}$.

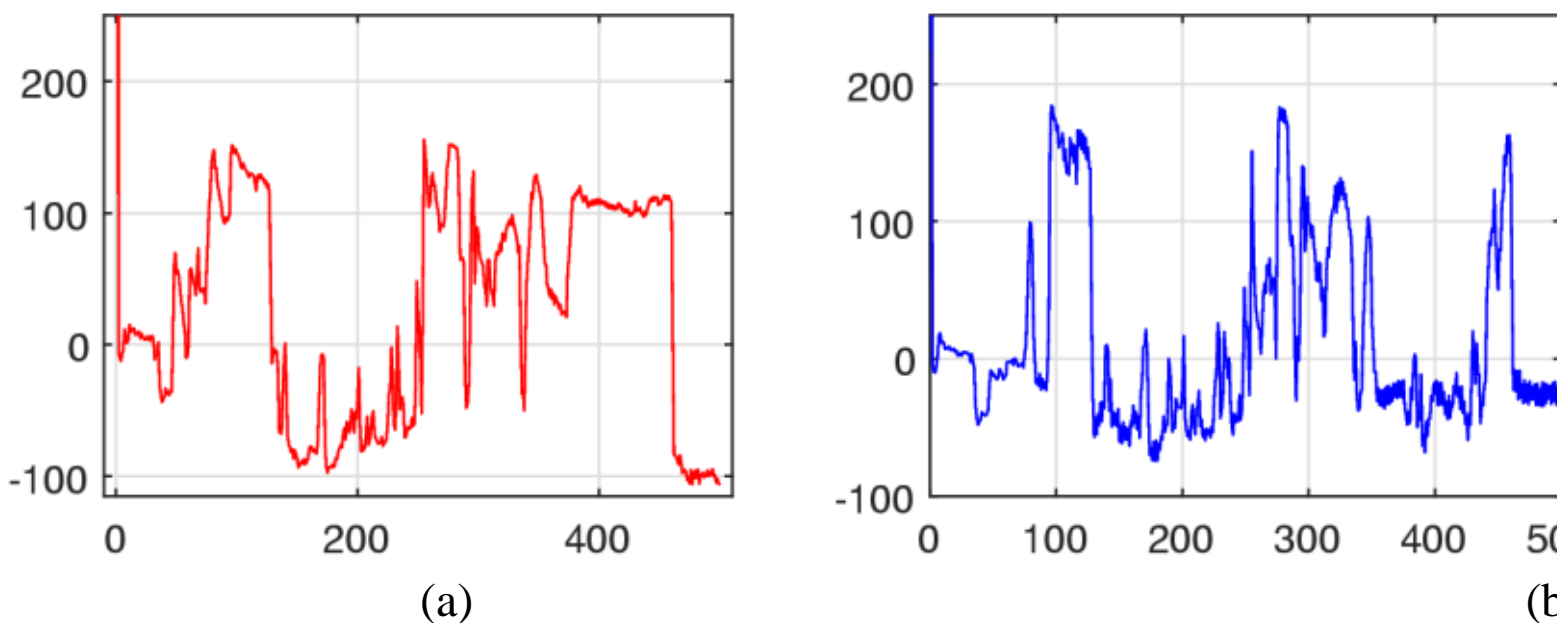

(a) (b)

**Figure 4:** The 500-point $\boldsymbol{T}$-type DsiHT of the signal $\boldsymbol{z}$; (a) real and (b) imaginary parts.

The results of the DsiHT of the signal $\boldsymbol{z}$, when using matrices of types $\boldsymbol{M}$ and $\boldsymbol{G}$ are almost the same as for the $\boldsymbol{T}$-type DsiHT. As mentioned above, the difference of these transforms is in the first two components;

$$\begin{aligned} \boldsymbol{T}-\text{based DsiHT} &= \{1962.3+389.2i, -5.3-3.9i, -9.8-9.6i, \dots .\}, \\ \boldsymbol{M}-\text{based DsiHT} &= \{1962.3+389.2i, -5.1-4.1i, -9.8-9.6i \ \dots .\}, \\ \boldsymbol{G}-\text{based DsiHT} &= \{1973.5+327.8i, -5.1-4.1i, -9.8-9.6i \ \dots .\}. \end{aligned}$$

Thus, we have three types of matrices **T**, **M**, and **G** for basic transforms that can be used in the QR decomposition by the DsiHT. These matrices are different, but these transforms move the energy of vector $\mathbf{x} = (\mathrm{x}_0, \mathrm{x}_1)$ to the first component and zeroing the second one.

## 4. QR DECOMPOSITIONS WITH THE DSIHT

In this section, we analyze the application of the DsiHTs that are based on basic transformations $\boldsymbol{T}$, $\boldsymbol{M}$ and $\boldsymbol{G}$ in QR decomposition of a square complex matrix. The QR-decomposition is described in detail on the example with a $4 \times 4$ matrix. Let $\boldsymbol{X}$ be the following $4 \times 4$ complex square matrix with $\det(\boldsymbol{X}) \neq 0$:

$$\boldsymbol{X} = \begin{bmatrix} a_0 & b_0 & c_0 & d_0 \\ a_1 & b_1 & c_1 & d_1 \\ a_2 & b_2 & c_2 & d_2 \\ a_3 & b_3 & c_3 & d_3 \end{bmatrix}.$$

First, we take the first column of the matrix as the vector $\boldsymbol{a} = (a_0, a_1, a_2, a_3)'$. This vecor will be the generator for the DsiHT, which we consider of the type $\boldsymbol{M}$. We denote the matrix of DsiHT by $\boldsymbol{T_a}$. The application of the DsiHT on the same vector is a vector $\overline{\boldsymbol{a}} = T_{\boldsymbol{a}}(\boldsymbol{a}) = (||\boldsymbol{a}||, 0, 0,\ 0)'$, where $||\boldsymbol{a}||$ is the energy of the vector, $||\boldsymbol{a}|| = \sqrt{|a_0|^2 + |a_1|^2 + |a_2|^2 + |a_3|^2}$. Therefore, when multiplying the matrices $\boldsymbol{T_a}$ and $\boldsymbol{X}$, we obtain the matrix $\boldsymbol{X}_1$ of the following form:

$$\boldsymbol{X}_1 = \boldsymbol{T_a}\boldsymbol{X} = \begin{bmatrix} ||\boldsymbol{a}|| & \hat{b}_0 & \hat{c}_0 & \hat{d}_0 \\ 0 & \hat{b}_1 & \hat{c}_1 & \hat{d}_1 \\ 0 & \hat{b}_2 & \hat{c}_2 & \hat{d}_2 \\ 0 & \hat{b}_3 & \hat{c}_3 & \hat{d}_3 \end{bmatrix}, \qquad \boldsymbol{Y}_1 = \begin{bmatrix} \hat{b}_1 & \hat{c}_1 & \hat{d}_1 \\ \hat{b}_2 & \hat{c}_2 & \hat{d}_2 \\ \hat{b}_3 & \hat{c}_3 & \hat{d}_3 \end{bmatrix}.$$

In the second step, this $3 \times 3$ complex submatrix $\boldsymbol{Y}_1$ is processed similarly. The first vector-column $\boldsymbol{b} = (\hat{b}_1, \hat{b}_2, \hat{b}_3)'$ is used as a generator for the 3-point DsiHT. We denote by $\boldsymbol{T_b}$ the matrix of this DsiHT. As a result, we obtain the following matrix:

$$\boldsymbol{X}_2 = \begin{bmatrix} 1 & 0 \\ 0 & \boldsymbol{T_b} \end{bmatrix} \boldsymbol{X}_1 = (1 \oplus \boldsymbol{T_b})\boldsymbol{X}_1 = (1 \oplus \boldsymbol{T_b})\boldsymbol{T_a}\boldsymbol{X} = \begin{bmatrix} ||\boldsymbol{a}|| & \hat{b}_0 & \hat{c}_0 & \hat{d}_0 \\ 0 & ||\boldsymbol{b}|| & \tilde{c}_1 & \tilde{d}_1 \\ 0 & 0 & \tilde{c}_2 & \tilde{d}_2 \\ 0 & 0 & \tilde{c}_3 & \tilde{d}_3 \end{bmatrix}$$

Here, $||\boldsymbol{b}|| = \sqrt{|\hat{b}_1|^2 + |\hat{b}_2|^2 + |\hat{b}_3|^2}$. In the last stage of calculation, the basic transform with the vector-generator $\boldsymbol{c} = (\tilde{c}_2, \tilde{c}_3)'$ is applied on two vectors $(\tilde{c}_2, \tilde{c}_3)'$ and $(\tilde{d}_2, \tilde{d}_3)'$. Denoting by $\boldsymbol{T_c}$ the matrix of this transform, we obtain the final triangularization,

$$\boldsymbol{R}=(1\oplus 1\oplus \boldsymbol{T}_c)\boldsymbol{X}_2=(1\oplus 1\oplus \boldsymbol{T}_c)(1\oplus \boldsymbol{T}_b)\boldsymbol{T}_a\boldsymbol{X}=\begin{bmatrix}||\boldsymbol{a}|| & \hat{b}_0 & \hat{c}_0 & \hat{d}_0\\ 0 & ||\boldsymbol{b}|| & \tilde{c}_1 & \tilde{d}_1\\ 0 & 0 & ||\boldsymbol{c}|| & \breve{d}_2\\ 0 & 0 & 0 & \breve{d}_3\end{bmatrix},$$

where $||\boldsymbol{c}||=\sqrt{|\tilde{c}_2|^2+|\tilde{c}_3|^2}$. The matrix of the transformation, or triangularization $\boldsymbol{X}\to\boldsymbol{R}$ is

$$\boldsymbol{T}=(1\oplus 1\oplus \boldsymbol{T}_c)(1\oplus \boldsymbol{T}_b)\boldsymbol{T}_{\boldsymbol{a}}. \tag{18}$$

Each of these DsiHTs is unitary, and therefore this matrix $\boldsymbol{T}$ is unitary. The inverse matrix $\boldsymbol{Q}=\boldsymbol{T}^{-1}$ is also unitary and can be written as

$$\boldsymbol{Q}=\boldsymbol{T}'=\boldsymbol{T}'_{\boldsymbol{a}}(1\oplus \boldsymbol{T}'_{\boldsymbol{b}})(1\oplus 1\oplus \boldsymbol{T}'_{\boldsymbol{c}}). \tag{19}$$

Thus, we have an explicit representation of the matrix $\boldsymbol{Q}$. Here, the operation $\boldsymbol{A}'$ denotes the conjugate transposition of a matrix $\boldsymbol{A}$**.** Thus, $\boldsymbol{TX}=\boldsymbol{R}$ and we obtain the following decomposition of the matrix $\boldsymbol{X}$:

$$\boldsymbol{X}=\boldsymbol{QR}. \tag{20}$$

It should be noted that if we apply instead of $\boldsymbol{T}$-type DsiHT the transform of type ***M*** or ***G*** in the above example, the diagonal coefficients of the matrix $\boldsymbol{R}$ will be changed as follows. In the case of the $\boldsymbol{M}$-type DsiHT,

$$||\boldsymbol{a}||\to \operatorname{sign}(a_0)||\boldsymbol{a}||,\qquad ||\boldsymbol{b}||\to \operatorname{sign}\left(\hat{b}_1\right)||\boldsymbol{b}||,\qquad ||\boldsymbol{c}||\to \operatorname{sign}(\tilde{c}_2)||\boldsymbol{c}||,$$

and in the case of the G-type DsiHT, these coefficients will be changed as

$$||\boldsymbol{a}||\to \operatorname{sign}(\boldsymbol{a})||\boldsymbol{a}||,\qquad ||\boldsymbol{b}||\to \operatorname{sign}(\boldsymbol{b})||\boldsymbol{b}||,\qquad ||\boldsymbol{c}||\to \operatorname{sign}(\boldsymbol{c})||\boldsymbol{c}||.$$

***Example 4×4:*** We consider the method of QR-decomposition that is similar to the method of the Householder transformations, only the DsiHTs will be used instead of the Householder transformations. First, we calculate the DsiHT by using the analytical equations, instead of matrix multiplications.

Let $\boldsymbol{X}$ be the following complex $4\times 4$ matrix:

$$\boldsymbol{X}=\begin{bmatrix}1+2i & 2-3i & 3+4i & -3+1i\\ 2-3i & 3+i & 2-2i & -6-7i\\ 1-i & 2-4i & 3+2i & 1+2i\\ 3-i & 4+3i & 4-2i & 2+4i\end{bmatrix}.$$

The method of QR decomposition of $\boldsymbol{X}$ with the DsiHTs results in the following presentations of the matrix $\boldsymbol{X}=\boldsymbol{QR}$.

*(a)* ***T-type DsiHT:***

The matrix $\boldsymbol{Q}$ is

$$\boldsymbol{Q}_T=\begin{bmatrix}0.1826+0.3651i & 0.3448-0.6035i & -0.2415+0.1577i & -0.5158+0.0299i\\ 0.3651-0.5477i & 0.0771+0.1906i & 0.0032-0.4489i & -0.5682+0.0075i\\ 0.1826-0.1826i & 0.1407-0.5490i & -0.0966-0.5316i & 0.5457-0.1495i\\ 0.5477-0.1826i & 0.2859+0.2677i & -0.4710+0.4489i & 0.2990+0.0224i\end{bmatrix}$$

and the triangular matrix equals

$$\boldsymbol{R}_T=\begin{bmatrix}5.4772 & 2.5560+2.7386i & 6.5727+0.5477i & 1.6432-1.4606i\\ 0 & 7.3462 & -1.6743+2.9403i & -2.7497+0.5763i\\ 0 & 0 & -3.3243 & 3.6995-4.9272i\\ 0 & 0 & 0 & 5.6893+6.0780i\end{bmatrix}.$$

The first column of the matrix $\boldsymbol{Q}_T$ is proportional to the first column of the matrix $\boldsymbol{X}$,

$$\begin{bmatrix} 0.1826+0.3651i \\ 0.3651-0.5477i \\ 0.1826-0.1826i \\ 0.5477-0.1826i \end{bmatrix} = 0.1826 \begin{bmatrix} 1+2i \\ 2-3i \\ 1-i \\ 3-i \end{bmatrix}.$$

***(b) M-type DsiHT***

The matrix $\boldsymbol{Q}$ is

$$\boldsymbol{Q}_M = \begin{bmatrix} 0.1826+0.3651i & 0.3448-0.6035i & 0.2415-0.1577i & -0.5166-0.0088i \\ 0.3651-0.5477i & 0.0771+0.1906i & -0.0032+0.4489i & -0.5671-0.0350i \\ 0.1826-0.1826i & 0.1407-0.5490i & 0.0966+0.5316i & 0.5554-0.1083i \\ 0.5477-0.1826i & 0.2859+0.2677i & 0.4710-0.4489i & 0.2999 \end{bmatrix}$$

and the triangular matrix is

$$\boldsymbol{R}_M = \begin{bmatrix} 5.4772 & 2.5560+2.7386i & 6.5727+0.5477i & 1.6432-1.4606i \\ 0 & 7.3462 & -1.6743+2.9403i & -2.7497+0.5763i \\ 0 & 0 & +3.3243 & -3.6995+4.9272i \\ 0 & 0 & 0 & 6.1279+5.6355i \end{bmatrix}.$$

Up to the signs $\pm$, many coefficients of the matrix $\boldsymbol{Q}_M$ are equal to the coefficients of $\boldsymbol{Q}_T$, and the main difference in the 4th columns of these matrices. In the matrices $\boldsymbol{R}_M$ and $\boldsymbol{R}_T$**,** the last coefficients are different, others are the same or differ only in the sign.

***(c) Householder Transforms***

We consider for comparison the Householder transform decomposition $\boldsymbol{X} = \boldsymbol{Q}_H \boldsymbol{R}_H$**.** The Householder transformation [6] is defined as $\boldsymbol{H} = \boldsymbol{I} - 2\boldsymbol{w}\boldsymbol{w}'$, where the normalized vector $\boldsymbol{w}$ is calculated

$$\boldsymbol{w} = \boldsymbol{x} - (\|\boldsymbol{x}\|, 0, 0, \dots, 0)' \;/\; \sqrt{2\|\boldsymbol{x}\|(\|\boldsymbol{x}\| - x_0)}$$

and $\boldsymbol{x} = (x_0, x_1, x_2, \dots, x_{N-1})'$ is a given vector. When running the MATLAB function '**qr.m**,' we obtain the following matrices:

$$\boldsymbol{Q}_H = \begin{bmatrix} -0.1826-0.3651i & -0.3448+0.6035i & 0.2415-0.1577i & 0.3743+0.3562i \\ -0.3651+0.5477i & -0.0771-0.1906i & -0.0032+0.4489i & 0.3937+0.4097i \\ -0.1826+0.1826i & -0.1407+0.5490i & 0.0966+0.5316i & -0.4821-0.2963i \\ -0.5477+0.1826i & -0.2859-0.2677i & 0.4710-0.4489i & -0.2207-0.2030i \end{bmatrix}$$

and

$$\boldsymbol{R}_H = \begin{bmatrix} -5.4772 & -2.5560-2.7386i & -6.5727-0.5477i & -1.6432+1.4606i \\ 0 & -7.3462 & 1.6743-2.9403i & 2.7497-0.5763i \\ 0 & 0 & 3.3243 & -3.6995+4.9272i \\ 0 & 0 & 0 & -8.3252 \end{bmatrix}.$$

The last column of the matrix $\boldsymbol{Q}_H$ is different from the last columns of the matrices $\boldsymbol{Q}_T$ and $\boldsymbol{Q}_M$, other coefficients are equal up to the sign. The same differences can be seen in the triangular matrices $\boldsymbol{R}_H$, $\boldsymbol{R}_T$, and $\boldsymbol{R}_M$.

***(d) G-type DsiHT:***

The matrix $\boldsymbol{Q}$ is

$$\boldsymbol{Q}_G = \begin{bmatrix} 0.4082 & 0.2457-0.6502i & -0.2362-0.1656i & -0.5166-0.0088i \\ -0.3266-0.5715i & 0.1061+0.1761i & 0.4179-0.1639i & -0.5671-0.0350i \\ -0.0816-0.2449i & 0.0527-0.5643i & 0.4576-0.2872i & 0.5554-0.1083i \\ 0.0816-0.5715i & 0.3244+0.2195i & -0.5918-0.2705i & 0.2999 \end{bmatrix}$$

and the triangular matrix is

$$\boldsymbol{R}_G = \begin{bmatrix} 2.4495+4.8990i & -1.3064+3.5109i & 2.4495+6.1237i & 2.0412+0.8165i \\ 0 & 7.2550+1.1542i & -2.1155+2.6407i & -2.8061+0.1371i \\ 0 & 0 & -1.2353+3.0863i & -3.1997-5.2656i \\ 0 & 0 & 0 & 6.1279+5.6355i \end{bmatrix}.$$

This decomposition is very different from the above three QR-decompositions. All coefficients of the first 3 columns of matrices $\boldsymbol{Q}_G$ and $\boldsymbol{R}_G$ differ from the coefficients of the corresponding matrices of the above decompositions by the $\boldsymbol{T}$-and $\boldsymbol{M}$-type, as well as the Householder QR decomposition. The last columns in matrices $\boldsymbol{Q}_G$ and $\boldsymbol{Q}_M$ are equal, as well as in matrices $\boldsymbol{R}_G$ and $\boldsymbol{R}_M$.

Here, we want to mention that if we apply the strong DsiHTs starting from the last column of the matrix, we obtain the decomposition of the matrix $\boldsymbol{X}=\boldsymbol{QL}$ with the left triangular matrix. For the $\boldsymbol{G}$-type strong DsiHT, the matrix $\boldsymbol{Q}$ is

$$\boldsymbol{Q}_G = \begin{bmatrix} 0.6434 & 0.1675-0.3605i & 0.5481-0.2101i & -0.0408+0.2858i \\ -0.1511-0.0403i & 0.1880+0.0742i & -0.1690-0.4448i & -0.8165+0.2041i \\ -0.6892+0.1466i & 0.2775-0.4503i & 0.3496-0.2445i & 0.2041 \\ 0.1270+0.2211i & 0.4693+0.5487i & -0.0886-0.4891i & 0.4082 \end{bmatrix}$$

and the left triangular matrix is

$$\boldsymbol{L}_G = \begin{bmatrix} -0.2137+1.5731i & 0 & 0 & 0 \\ 1.1871-1.9594i & 7.9344-0.8122i & 0 & 0 \\ 1.9415+4.1538i & 0.6302+0.7221i & 2.5389+7.6166i & 0 \\ -0.2858+1.0614i & -1.1431-1.4697i & 1.2247-0.2041i & 4.8990+9.7980i \end{bmatrix}.$$

One can notice from this example that the QR-decomposition by the $\boldsymbol{G}$-type DsiHTs differs much from the decomposition by the Householder transforms. The QR-decomposition by $\boldsymbol{M}$-type DsiHTs is very similar to the Householder transform method. The $\boldsymbol{M}$-type DsiHT is fast and can be implemented by using 2×2 basic transformations, as well by the analytical equations. To analyse the difference of these two decompositions in more detail, we consider the example with a 6×6 complex matrix $\boldsymbol{X}$.

***Example 6×6:*** Let the matrix $\boldsymbol{X}$ be the following:

$$\boldsymbol{X} = \begin{bmatrix} 1+2i & 2-3i & 3+4i & -3+i & -4-i & 2-3i \\ 2-3i & 3+i & 2-2i & -6-7i & 2+i & 5-2i \\ 4-i & 3-2i & 4-5i & 2+3i & 4+7i & 6+2i \\ 5+2i & 5+i & 3-2i & 8-3i & 7-2i & 2+3i \\ 4-3i & -5-2i & 1-i & 2-4i & 3+2i & 1+2i \\ 7-2i & 6+i & 3-i & 4+3i & 4+3i & 2+4i \end{bmatrix}.$$

The QR-decomposition by the $\boldsymbol{M}$-type DsiHTs results in the following matrices:

$$\boldsymbol{Q}_M = \begin{bmatrix} 0.0839 & 0.1419 & 0.3046 & -0.1235 & 0.4129 & -0.0136 \\ 0.1678 & 0.2321 & 0.2051 & -0.4530 & 0.2174 & 0.4402 \\ 0.3357 & 0.1232 & 0.2883 & -0.0854 & 0.0096 & -0.4182 \\ 0.4196 & 0.2579 & -0.0885 & 0.3133 & 0.3378 & -0.1404 \\ 0.3357 & -0.6763 & -0.0301 & -0.0356 & 0.3404 & -0.2840 \\ 0.5874 & 0.2937 & -0.1343 & 0.0922 & -0.1331 & 0.1813 \end{bmatrix} +$$

$$+i\begin{bmatrix} 0.1678 & -0.3926 & 0.5185 & -0.0103 & -0.2061 & 0.4476 \\ -0.2518 & 0.2579 & 0.0226 & -0.4400 & 0.2949 & -0.1367 \\ -0.0839 & -0.1275 & -0.5164 & 0.1108 & 0.4122 & 0.3669 \\ 0.1678 & 0.0430 & -0.2965 & -0.3644 & -0.4614 & -0.2323 \\ -0.2518 & -0.0330 & 0.2352 & 0.0055 & 0.1455 & -0.3004 \\ -0.1678 & 0.2465 & 0.2758 & 0.5705 & -0.0318 & 0 \end{bmatrix}$$

and

$$
\boldsymbol{R}_M = \begin{bmatrix} 11.9164 & 5.5386 & 6.9652 & 7.4687 & 6.1260 & 4.5316 \\ 0 & 9.8295 & 0.6133 & -1.5246 & 0.4542 & 4.0665 \\ 0 & 0 & 6.4709 & -3.2862 & -4.5013 & 0.9395 \\ 0 & 0 & 0 & 11.9062 & 1.6459 & 0.0832 \\ 0 & 0 & 0 & 0 & 6.3390 & 2.3524 \\ 0 & 0 & 0 & 0 & 0 & -2.1708 \end{bmatrix} +
$$

$$
+\, i \begin{bmatrix} 0 & -0.8392 & -2.8532 & -1.9301 & 2.8532 & 6.0421 \\ 0 & 0 & -0.3095 & 1.0603 & -4.2671 & -0.3324 \\ 0 & 0 & 0 & 3.2384 & 6.6643 & 0.1439 \\ 0 & 0 & 0 & 0 & -1.1619 & 3.4811 \\ 0 & 0 & 0 & 0 & 0 & -3.1871 \\ 0 & 0 & 0 & 0 & 0 & -3.5886 \end{bmatrix}.
$$

The method of Householder transforms which is performed by the MATLAB function ‘**`qr.m`**’ results in the following matrices $\boldsymbol{Q}_H$ and $\boldsymbol{R}_H$ in the QR-decomposition $\boldsymbol{X} = \boldsymbol{Q}_H\boldsymbol{R}_H$. In the matrix $\boldsymbol{Q}_H$, the last column differs from the same column in the matrix $\boldsymbol{Q}_M$, when using the $\boldsymbol{M}$-type matrices in QR-decomposition. The remaining part, or the 6×5 sub-matrix of $\boldsymbol{Q}_H$ differs from the same sub-matrix of $\boldsymbol{Q}_M$ only by the sign. Thus, we can write the matrix $\boldsymbol{Q}_H$ in the following way:

$$
\boldsymbol{Q}_H = -(\boldsymbol{Q}_M)_{6\times 5} \sqcup \underbrace{\begin{bmatrix} -0.3900 \\ 0.3448 \\ -0.5304 \\ 0.1261 \\ 0.1100 \\ 0.0939 \end{bmatrix} + i \begin{bmatrix} 0.2201 \\ 0.3059 \\ -0.1679 \\ -0.2404 \\ -0.3985 \\ 0.1552 \end{bmatrix}}.
$$

All coefficients of the triangular matrix $\boldsymbol{R}_H$, except the last coefficient, equal to the corresponding coefficients of the matrix $\boldsymbol{R}_M$ with the sign minus. In matrix $\boldsymbol{R}_H$, this coefficient equals $-4.1941$ and in the matrix $\boldsymbol{R}_M$ such coefficient is $2.1708 + 3.5886$. Therefore, we can write that

$$
\boldsymbol{R}_H = -\boldsymbol{R}_M + \begin{bmatrix} 0 \\ 0 \\ 0 \\ 0 \\ 0 \\ 4.1941 - (2.1708 + 3.5886) \end{bmatrix}.
$$

## 5. THE SCRIPTS FOR THE DECOMPOSITION BY THE DSIHT

Below are MATLAB-based codes for QR-decomposition of a complex square matrix $\boldsymbol{X}$, by using the DsiHT. The function ‘**`amqr_compex.m`**’ is for QR-decomposition of $\boldsymbol{X}$ by the DsiHT with the basic transforms $\boldsymbol{T}$, $\boldsymbol{M}$, and $\boldsymbol{G}$. The Householder transform-based QR decomposition is calculated in this example by the MATLAB function ‘**`qr.`**’ All matrices in the examples given in Section 4 were calculated by these codes.

```
% ---------------------------------------------------------------------------
% call: run_Matrix4x4general.m / from the library of Art Grigoryan
% QR-decomposition by the DsiHT (of 3 types)
  X=[1+2j  2-3j 3+4j -3+1j
     2-3j  3+1j 2-2j -6-7j
     1-1j  2-4j 3+2j  1+2j
     3-1j  4+3j 4-2j  2+4j];
  ntype=1;
  [U,R]=amqr_complex(X,ntype);  % the DsiHT method of QR-decomposition
  Q=U';
  Q*R % = X
  [QH,RH]=qr(X);     % MATLAB function (using the Householder transform)
  QH*RH % = X
% ---------------------------------------------------------------------------
```

The function **amqr_complex** is used to calculate the QR decomposition by the DsiHTs. The parameter '**ntype**' is selected as 1, 2, and 3 for the ***T***, ***M***, and ***G*** -type DsiHTs, respectively. If parameter '**ntype**' is set to 0, the ***M*** -type DsiHT is calculated by analytical equations (4)-(6) in the function **msob_enanalcomp**. The calculation of 2×2 basic matrices of types ***T***, ***M***, and ***G*** are accomplished by the corresponding functions **msob2T**, **msob2M**, and **msob2G** in the main function **msob_complexA** that calculates the DsiHT. The scripts of all these functions are given below.

```
% -----------------------------------------------------------------
% call: amqr_complex / from the library of Art Grigoryan, 2014-2019
% The QR decomposition of NxN matrix A by the DsiHT of type:
%  #1, for T-type DsiHT / using 2x2 matrices T by msob2T.m
%  #0, for M-type DsiHT / analytical by msob_enanalcomp.m
%  #2, for M-type DsiHT / using 2x2 matrices M by msob2M.m
%  #3, for G-type DsiHT / using 2x2 matrices G by msob2G.m
function [U,A]=amqr_complex(A,ntype)
  N=size(A,1);
  E=eye(N);  U=E;
  for j1=1:N-1
      E=eye(N);
      A1=A(j1:N,j1);
      for i1=j1:N
          if (ntype==0)
              y(j1:N,i1)=msob_enanalcomp(A1,A(j1:N,i1));
              E(j1:N,i1)=msob_enanalcomp(A1,E(j1:N,i1));
          else
              y(j1:N,i1)=msob_complexA(A1,A(j1:N,i1),ntype);
              E(j1:N,i1)=msob_complexA(A1,E(j1:N,i1),ntype);
          end
      end
      A=y;
      U=E*U;
  end
end

% -----------------------------------------------------------------
% call: msob_analcomp.m / from the library of Art Grigoryan
% x-induced complex discrete heap transform, z=T_x[y].
% x is the complex generator, and y is the input signal.
  function z=msob_enanalcomp(x,y)
    N=length(x);
    Ex=x(1);  Exc=conj(Ex);
    Ey=y(1);  Eyc=conj(Ex);
    Eyx=Ey*Exc;
    Ex2=Ex*Exc;
    for i1=2:N
        Ex=sqrt(Ex2);
        xi1=x(i1);  xi1c=conj(xi1);
        yi1=y(i1);  yi1c=conj(yi1);
        zi1=(Eyx/Ex)*xi1 - yi1*Ex;
        Eyx=Eyx+yi1*xi1c;
        Ex2=Ex2+xi1*xi1c;
        Ex=sqrt(Ex2);
        zi1=-zi1/Ex;
        z(i1)=zi1;
    end
    z(1)=Eyx/Ex;
  end

% -------------------------------------------------------------------
% call: msob_complexA.m / from the library of Art Grigoryan
% x-vector generated complex Discrete Heap Transform of the vector Q.
% Transform type: T (#1), M(#2), and G(#3).  x-DsiHT: Q -> Z
  function Z=msob_complexA(X,Q,ntype)
    N=size(X,1);
```

```
    x1=X(1);
    q1=Q(1);
    Z=zeros(N,1);
    for i1=2:N
        x2=X(i1);
        q2=Q(i1);
        if (ntype==1)
            [d1,d2]=msob2T(x1,x2,q1,q2);
            x1=sign(real(x1))*norm([x1,x2]);
        elseif (ntype==2)
            [d1,d2]=msob2M(x1,x2,q1,q2);
            x1=norm([x1,x2]);
        else
            [d1,d2]=msob2G(x1,x2,q1,q2);
            x1=x1*sqrt(1 + norm(x2)^2 / norm(x1)^2 );
        end
        q1=d1;
        Z(i1)=d2;
    end
    Z(1)=q1;
 end

% -----------------------------------------------------
% call: msob2T.m   / from the library of Art Grigoryan
% to calculate 2x2 basic transform with complex matrix T
% generator: 2 complex vectors x,y
% input / output: 2 complex vectors c1,c2 / d1,d2
function [d1,d2]=msob2T(x,y,c1,c2)
  T = [ conj(x)   conj(y)
          -y        x     ];
  T = sign(real(x))*T;
  T=T/sqrt( norm(x)^2 + norm(y)^2 );
  z=[c1; c2];
  d=T*z;
  d1=d(1); d2=d(2);

% -----------------------------------------------------
% call: msob2M.m  / from the library of Art Grigoryan
% to calculate 2x2 basic transform with complex matrix M
function [d1,d2]=msob2M(x,y,c1,c2)
  x_conj=conj(x);
  y_conj=conj(y);
  d1 = x_conj*c1 + y_conj*c2;
  nx=norm(x);
  yx = y*x_conj/nx;
  d2 = -yx*c1 + nx*c2;
  E20=sqrt( nx^2 +norm(y)^2 );    % for generator
  d1=d1/E20; d2=d2/E20;
end

% -------------------------------------------------------
% call: msob2G.m  / from the library of Art Grigoryan
% to calculate 2x2 basic transform with complex matrix G
function [d1,d2]=msob2G(x,y,c1,c2)
  E20=sqrt(norm(x)^2 +norm(y)^2 );
  xy = x*conj(y);
  nx=norm(x);
  d1 = nx*c1 + xy*c2/nx;
  yx = y*conj(x);
  d2 = -yx*c1/nx + nx*c2;
  d1=d1/E20; d2=d2/E20;
end
% -------------------------------------------------------
```

Below is the script of the function '**msob_complex1sp.m**' to compute the strong $\boldsymbol{G}$-type DsiHT. To calculate the strong $\boldsymbol{T}$ and $\boldsymbol{M}$-types DsiHT, the corresponding function can be written in a similar way.

```
% -------------------------------------------------------------
% call: msob_complex1sp.m   / from the library of Art Grigoryan
% To calculate the strong G-type X-generated DsiHT: Q->Z.
function Z=msob_complex1sp(X,Q)
    N=size(X,1);
    x1=X(N);
    q1=Q(N);
    Z=zeros(N,1);
    for i1=N-1:-1:1
        x2=X(i1);
        q2=Q(i1);
        [d1,d2]=msob2G(x1,x2,q1,q2);  % can use  msob2T, msob2M
        q1=d1;
        Z(i1,:)=d2;
        x1=x1*sqrt(1 + norm(x2)^2 / norm(x1)^2 );  % !
    end
    Z(N)=q1;
end
```

## 6. Mixed Type DsiHT QR Decomposition

It is clear that on different stages of QR-decomposition by the DsiHT, we can change the type of the DsiHT. Such mixed type QR-decompositions can be considered and applied in practice together with the described above QR-decompositions by the $\boldsymbol{T}$, $\boldsymbol{M}$, and $\boldsymbol{G}$-type DsiHTs.

As an example, we consider the matrix $\boldsymbol{X}$ given in Example 6×6 in Section 4 with the following set of types of transforms: [T M G T T], or [1 2 3 1 1] in numbers when running the codes. It means that the first transform which will be used to obtain the first heap in the first column of in the matrix $\boldsymbol{X}$ is the $\boldsymbol{T}$-type DsiHT. The second transform is the $\boldsymbol{M}$-type DsiHT to get the second heap in the matrix, or coefficient number (2,2) in the triangular matrix $\boldsymbol{R}$, and so on. As a result, we obtain the decomposition of the matrix $\boldsymbol{X}=\boldsymbol{Q}_{12311}\boldsymbol{R}_{12311}$, where the unitary matrix is

$$\boldsymbol{Q}_{12311}=\begin{bmatrix} 0.0839 & 0.1419 & -0.5953 & -0.1235 & -0.4129 & 0.3665 \\ 0.1678 & 0.2321 & -0.0986 & -0.4530 & -0.2174 & 0.1278 \\ 0.3357 & 0.1232 & 0.3685 & -0.0854 & -0.0096 & 0.0767 \\ 0.4196 & 0.2579 & 0.3079 & 0.3133 & -0.3378 & -0.2713 \\ 0.3357 & -0.6763 & -0.2063 & -0.0356 & -0.3404 & -0.4070 \\ 0.5874 & 0.2937 & -0.2043 & 0.0922 & 0.1331 & 0.0997 \end{bmatrix}+$$

$$+i\begin{bmatrix} 0.1678 & -0.3926 & 0.0853 & -0.0103 & 0.2061 & 0.2573 \\ -0.2518 & 0.2579 & 0.1812 & -0.4400 & -0.2949 & -0.4429 \\ -0.0839 & -0.1275 & 0.4625 & 0.1108 & -0.4122 & 0.5510 \\ 0.1678 & 0.0430 & 0.0305 & -0.3644 & 0.4614 & -0.0103 \\ -0.2518 & -0.0330 & -0.1170 & 0.0055 & -0.1455 & 0.0722 \\ -0.1678 & 0.2465 & -0.2289 & 0.5705 & 0.0318 & -0.1515 \end{bmatrix}.$$

When comparing with the QR-decomposition by the $\boldsymbol{M}$-type DsiHT, one can notice that the columns number 3 and 6 and the sign on the 5th column are different in matrices $\boldsymbol{Q}_{12311}$ and $\boldsymbol{Q}_{\boldsymbol{M}}$. The triangular matrix equals to

$$\boldsymbol{R}_{12311}=\begin{bmatrix} 11.9164 & 5.5386 & 6.9652 & 7.4687 & 6.1260 & 4.5316 \\ 0 & 9.8295 & 0.6133 & -1.5246 & 0.4542 & 4.0665 \\ 0 & 0 & -2.4534 & 4.2425 & 7.8733 & -0.2230 \\ 0 & 0 & 0 & 11.9062 & 1.6459 & 0.0832 \\ 0 & 0 & 0 & 0 & 6.3390 & -2.3524 \\ 0 & 0 & 0 & 0 & 0 & 1.8050 \end{bmatrix}+$$

$$+\,i\begin{bmatrix} 0 & -0.8392 & -2.8532 & -1.9301 & 2.8532 & 6.0421 \\ 0 & 0 & -0.3095 & 1.0603 & -4.2671 & -0.3324 \\ 0 & 0 & -5.9878 & 1.8131 & 1.6386 & -0.9239 \\ 0 & 0 & 0 & 0 & -1.1619 & 3.4811 \\ 0 & 0 & 0 & 0 & 0 & 3.1871 \\ 0 & 0 & 0 & 0 & 0 & -3.7858 \end{bmatrix}.$$

In this matrix, the coefficients of the 3rd and 6th rows and the sign in 5th row differ from the corresponding coefficients of matrix $\boldsymbol{R_M}$.

It should be noted that instead of combination of types [1 2 3 1 1], we can use other combinations with 1,2, and 3. The number of such combinations is $3^5 = 243$. In general case of the $N \times N$ matrix, we can choose one combination of with numbers presenting the types of the DsiHT. The number of such combinations equals $3^{(N-1)}$ and they can be used for calculating the QR-decomposition by the DsiHTs. The combinations with all 1s, 2s, and 3s correspond to the QR-decompositions by the $\boldsymbol{T}$, $\boldsymbol{G}$, and $\boldsymbol{M}$-type DsiHTs. Also, different paths can be used for the DsiHT, which increases the number of possible QR decompositions.

In conclusion, we consider a few QR-decompositions which were calculated for the pseudorandom integer $N \times N$ matrices $\boldsymbol{X}$ with complex coefficients with real and imaginary parts in the range $1{:}\,N$. For that, the MATLAB function '**randi.m**' is used. Twelve values of $N$ were arbitrary chosen to be 6, 13, 17, 19, 21, 40, 64, 100, 128, 201, 256, and 400. The QR-decompositions for each of these values of $N$ were calculated by the $\boldsymbol{M}$-type DsiHT and the Householder transforms; the script of the code is given below for the QR decomposition of the random complex 400× 400 matrix.

```
% ======================================================
% call: testQR.com    / Art Grigoryan
  N=400;
  X=randi(N,N)+randi(N,N)*j;
  ntype=0;
  [U,R]=amqr_complex(X,ntype);
  Q=U'; Y=Q*R;    % = X
  E1=norm(X-Y)    % 4.8725e-11
  [Qm,Rm]=qr(X); % MATLAB: QR by the Householder transf.
  Ym=Qm*Rm;       % = X
  E2=norm(X-Ym)   % 5.5772e-11
% ======================================================
```

# 6. Summary Results

To compare the results of the QR-decomposition, the precision of computation was estimated by the 2-norms of the matrix $\Delta_{\boldsymbol{M}} = (\boldsymbol{X} - \boldsymbol{Q_M R_M})$ and matrix $\Delta_{\boldsymbol{H}} = (\boldsymbol{X} - \boldsymbol{Q_H R_H})$, by using the MATLAB function "**norm.m.**" The results of estimations are given in Table 1.

| $N$ | norm ($\Delta_{\boldsymbol{M}}$) | norm ($\Delta_{\boldsymbol{H}}$) | $N$ | norm ($\Delta_{\boldsymbol{M}}$) | norm ($\Delta_{\boldsymbol{H}}$) |
|---|---|---|---|---|---|
| 6 | 5.0854e-15 | 1.5907e-14 | 64 | 7.9044e-13 | 1.3359e-12 |
| 13 | 2.8659e-14 | 5.9598e-14 | 100 | 2.4268e-12 | 2.9176e-12 |
| 17 | 4.8721e-14 | 6.2739e-14 | 128 | 4.8050e-12 | **4.7850e-12** |
| 19 | 5.7744e-14 | 7.0202e-14 | 201 | 1.0487e-11 | 2.1172e-11 |
| 21 | 9.2941e-14 | 1.2370e-13 | 256 | 1.6789e-11 | 2.8667e-11 |
| 40 | 3.6162e-13 | 4.7966e-13 | 400 | 4.8725e-11 | 5.5772e-11 |

Table 1. The precision of the QR-decomposition by the $\boldsymbol{M}$-type DsiHT and Householder transforms.

One can notice that in most cases the 2-norm of the QR-decomposition by the $\boldsymbol{M}$-type DsiHTs is less than the same norm when using the Householder transforms.

## Conclusion

We described three different types of QR-decompositions, that includes the DsiHT with $\boldsymbol{T}$, $\boldsymbol{G}$, and $\boldsymbol{M}$-type complex matrices. The decomposition by analytical formulas was also given for the $\boldsymbol{M}$-type DsiHT. The mixed type QR-decomposition, when different type DsiHTs are used in different stages of the algorithm was also presented. For an $N \times N$ complex nonsingular matrix, the number of such decompositions is greater than $3^{(N-1)}$. Examples of the QR-decomposition were given in detail for the 4×4 and 6×6 complex matrices and compared with the Householder transform-based QR-decomposition. MATLAB-based scripts of the codes for calculating the DsiHTs and QR decompositions are given. The different QL-decompositions of a complex matrix can be obtained in a similar way. We believe that the concept of the DsiHT can be generalized to the quaternion space [20] and the QR-decomposition of the quaternion matrix can be calculated and used together with the known methods of decompositions [19].

## References

[1] D.D. Morrison, *Remarks on the Unitary Triangularization of a Nonsymmetric Matrix*. Journal ACM, **7** (2) (1960) 185-186.
[2] P. Businger, G.H. Golub, *Linear Least Squares Solutions by Householder Transformations*. In: Bauer F.L. (eds) Linear Algebra. Handbook for Automatic Computation, **2** (2017). Springer, Berlin, Heidelberg.
[3] R.A. Horn, R.J. Charles, *Matrix Analysis*. Cambridge University Press, 1985.
[4] G.H. Golub, C.F. Van Loan, *Matrix Computations*, 3rd edition, Johns Hopkins, 1996.
[5] E. Anderson, Z. Bai, C.H. Bischof, et al., *LAPACK Users' Guide*. 3rd edition, SIAM Philadelphia, PA, 1999.
[6] T.K. Moon, W.C. Stirling, *Mathematical Methods and Algorithms for Signal Processing*. NJ: Prentice Hall, 2000.
[7] P. Alonso, J.M. Peña, M.L. Serrano, *QR decomposition of almost strictly sign regular matrices*. Journal of Computational and Applied Mathematics, **318** (2017), 646-657.
[8] A. Björck, *Solving Linear Least Squares Problems by Gram-Schmidt Orthogonalization*, BIT **7** (1967) 1-21.
[9] H. Rutishause, *Simultaneous Iteration Method for Symmetric Matrices*. In: In: Bauer F.L. (eds) Linear Algebra. Handbook for Automatic Computation, **186** (2017). Springer, Berlin, Heidelberg.
[10] A.S. Householder, Unitary Triangulation of a Nonsymmetric Matrix. Journal ACM, **5** (4), (1958) 339-342.
[11] D. Bindel, J. Demmel, W. Kahan, O. Marques, *On Computing Givens rotations reliably and efficiently*. LAPACK Working Note 148, University of Tennessee, UT-CS-00-449, 2001.
[12] J. Demmel, L. Grigori, M. Hoemmen, J. Langou, *Communication-optimal parallel and sequential QR and LU factorizations*, SIAM J. Sci. Comp., **34** (1), (2012) 206-239.
[13] A.M. Grigoryan, M.M. Grigoryan, *Nonlinear approach of construction of fast unitary transforms*. Proceedings of the 40th Annual Conference on Information Sciences and Systems (CISS 2006), Princeton University, Princeton, (2006) 1073-1078.
[14] A.M. Grigoryan, M.M. Grigoryan, *Discrete unitary transforms generated by moving waves*. Proceedings of the International Conference: Wavelets XII, SPIE: Optics + Photonics 2007, San Diego, CA, 2007.
[15] A.M. Grigoryan, M.M. Grigoryan, *New discrete unitary Haar-type heap transforms*. Proceedings of the International Conference: Wavelets XII, SPIE: Optics + Photonics 2007, San Diego, CA, 2007.
[16] A.M. Grigoryan, M.M. Grigoryan, *Brief Notes in Advanced DSP: Fourier Analysis with MATLAB*. CRC Press, Taylor and Francis Group, 2009.
[17] A.M. Grigoryan, *New method of Givens rotations for triangularization of square matrices*, Journal of Advances in Linear Algebra & Matrix Theory (ALAMT), **4** (2), (2014) 65-78.
[18] A.M. Grigoryan, *Fast Heap transform-based QR-decomposition of real and complex matrices: Algorithms and codes*, [9411-21], Proceedings of SPIE vol. 9411, Electronic Imaging: Mobile Devices and Multimedia: Enabling Technologies, Algorithms, and Applications 2015, February 10-11, San Francisco, California.
[19] L. Ying, W. Musheng, Z. Fengxia, Z. Jianli, *Real structure-preserving algorithms of Householder based transformations for quaternion matrices*, Journal of Computational and Applied Mathematics, **305** (2016), 82-91.
[20] A.M. Grigoryan and S.S. Agaian, *Quaternion and Octonion Color Image Processing with MATLAB*, p. 404, SPIE, vol. PM279, April 5, 2018.